\documentclass[a4paper,11pt]{article}

\usepackage{amsfonts,amssymb,amsbsy,amsmath,amsthm,mathrsfs,enumerate,verbatim}

  \usepackage{graphics} 
  \usepackage{epsfig} 
\usepackage{graphicx}  \usepackage{epstopdf}
 \usepackage[colorlinks=true]{hyperref}

\hypersetup{urlcolor=blue, citecolor=red}

\usepackage{amsfonts}
\usepackage{amsmath}
\usepackage{amssymb}
\usepackage{graphicx}
\usepackage{epsf}
\usepackage{float}

\usepackage{epsfig}

\usepackage{color}
\usepackage{afterpage}
\usepackage{verbatim}

\usepackage{times}

\newtheorem{theo}{Theorem}

\newtheorem{theorem}{Theorem}[section]

\theoremstyle{definition}
\newtheorem{remark}[theorem]{Remark}

\newcommand{\bel}{\begin{equation} \label}
\newcommand{\ee}{\end{equation}}

\def\beq{\begin{equation}}
\def\eeq{\end{equation}}
\newcommand{\bea}{\begin{eqnarray}}
\newcommand{\eea}{\end{eqnarray}}
\newcommand{\beas}{\begin{eqnarray*}}
\newcommand{\eeas}{\end{eqnarray*}}

{

 \usepackage{graphicx,color}
 \definecolor{mygreen}{cmyk}{1,0,1,0.1}

\usepackage{epsf}
\usepackage{epstopdf}

\usepackage{pdfsync}

\graphicspath{
}

\begin{document}

\title{Numerical validation of an explicit $P_1$ finite-element scheme for
  Maxwell's equations in a polygon with variable permittivity away from its boundary}

\author{L. Beilina  \thanks{
Department of Mathematical Sciences, Chalmers University of Technology and
 University of Gothenburg, SE-42196 Gothenburg, Sweden, e-mail: \texttt{\
larisa@chalmers.se}}
\and
V. Ruas \thanks{Institut Jean Le Rond d'Alembert, UMR 7190 CNRS - Sorbonne Universit\'e, F-75005, Paris, France, e-mail: \texttt{\
vitoriano.ruas@upmc.fr}}
}

\date{}

\maketitle

\begin{abstract}

  This paper is devoted to the numerical validation of an explicit
  finite-difference scheme for the integration in time of Maxwell's
  equations in terms of the sole electric field, using standard linear
  finite elements for the space discretization. The rigorous
  reliability analysis of this numerical model was the object of
  another authors' arXiv paper. More specifically such a study applies
  to the particular case where the electric permittivity has a
  constant value outside a sub-domain, whose closure does not
  intersect the boundary of the domain where the problem is
  defined. Our numerical experiments in two-dimension space certify
  that the convergence results previously
  derived for this approach are optimal, as long as the underlying CFL
  condition is satisfied.
\end{abstract}

\maketitle


\section{Motivation}

The purpose of this article is to provide a numerical validation of an
explicit $P_1$ finite element solution scheme of hyperbolic Maxwell's
equations for the electric field with constant dielectric permittivity
in a neighborhood of the boundary of the computational domain.  This
numerical model was thoroughly studied in \cite{arxiv1} from the
theoretical point of view. The focus of that paper was on the case of
three-dimensional domains on whose boundary absorbing boundary
conditions are enforced. However as pointed out therein, all the
underlying analytical results trivially extend to the two-dimensional
case and to other types of boundary conditions such as Dirichlet and
Neumann conditions. Actually in our experiments we consider only
two-dimensional test-problems, in which Dirichlet boundary conditions
are prescribed.  \\

The standard continuous $P_1$ FEM is a tempting possibility to solve
Maxwell's equations, owing to its simplicity.  It is well known
however that, for different reasons, this method is not always well
suited for this purpose. The first reason is that in general the
natural function space for the electric field is not the Sobolev space
${\bf H}^1$, but rather in the space ${\bf H}(curl)$. Another issue
difficult to overcome with continuous Lagrange finite elements is the
prescription of the zero tangential-component boundary conditions for
the electric field, which hold in many important applications. All
this motivated the proposal by N\'ed\'elec about four decades ago of a
family of ${\bf H}(curl)$-conforming methods to solve these equations
(cf. \cite{Nedelec}). These methods are still widely in use, as much
as other approaches well adapted to such specific conditions (see
e.g. \cite{Assous}, \cite{CiarletJr} and \cite{AMIS}). A comprehensive
description of finite element methods for Maxwell's equations can be
found in \cite{Monk}. \\ \indent There are situations however in which
the $P_1$ finite element method does provide an inexpensive and
reliable way to solve the Maxwell's equations. In this work we address
one of such cases, characterized by the fact that the electric
permittivity is constant in a neighborhood of the whole boundary of
the domain of interest. This is because, at least in theory, whenever
the electric permittivity is constant, the Maxwell's equations
simplify into as many wave equations as the space dimension under
consideration.
More precisely here we show by means of numerical
examples that, in such a particular case, a space discretization with
conforming linear elements, combined with a straightforward explicit
finite-difference scheme for the time integration, gives rise to
optimal approximations of the electric field, as long as a classical
CFL condition is satisfied.\\

\indent Actually this work is strongly
connected with studies presented in \cite{BG, cejm} for a combination
of the finite difference method in a sub-domain with constant
permittivity with the finite element method in the complementary
sub-domain. As pointed out above, the Maxwell's equations reduces to
the wave equation in the former case. Since the analysis of
finite-difference methods for this type of equation is well
established, only an explicit $P_1$ finite element scheme for
Maxwell's equations is analyzed in this paper. \\ \indent In \cite{BG,
  cejm} a stabilized domain-decomposition
finite-element/finite-difference approach for the solution of the
time-dependent Maxwell's system for the electric field was proposed
and numerically verified.  In these works \cite{BG, cejm} different
manners to handle a divergence-free condition in the finite-element
scheme were considered. The main idea behind the domain decomposition
methods in \cite{BG, cejm} is that a rectangular computational domain
is decomposed into two sub-domains, in which two different types of
discretizations are employed, namely, the finite-element domain in
which a classical $P_1$ finite element discretization is used, and the
finite-difference domain, in which the standard five- or seven-point
finite difference scheme is applied, according to the space
dimension. The finite element domain lies strictly inside the finite
difference domain, in such a way that both domains overlap in two
layers of structured nodes. First order absorbing boundary conditions
\cite{EM} are enforced on the boundary of the computational domain,
i.e. on the outer boundary of the finite-difference domain. In
\cite{BG, cejm} it was assumed that the dielectric permittivity
function is strictly positive and has a constant value in the
overlapping nodes as well as in a neighborhood of the boundary of the
domain. An explicit scheme was used both in the finite-element and
finite-difference domains. \\
\indent We recall that for a stable finite-element
solution of Maxwell's equation divergence-free edge elements are the
most satisfactory from a theoretical point of view ~\cite{Nedelec,
  Monk}.  However, the edge elements are less attractive for solving 
time-dependent problems, since a linear system of equations should
be solved at every time iteration.  In contrast, $P_1$ elements can be
efficiently used in a fully explicit finite element scheme with lumped
mass matrix \cite{delta, joly}. On the other hand it is also well known that the 
numerical solution of Maxwell's equations with nodal finite elements can 
result in unstable spurious solutions \cite{MP, PL}. Nevertheless a number of techniques 
are available to remove them, and in this respect we refer for
example to \cite{Jiang1, Jiang2, Jin, div_cor, PL}. 
In the current work, similarly to \cite{BG, cejm}, the spurious
solutions are removed  from the finite element scheme by adding the
divergence-free condition to the model equation for the electric
field.  Numerical tests given in \cite{cejm} demonstrate that 
spurious solutions are removable indeed, in case an explicit $P_1$ finite-element
solution scheme is employed.\\
\indent Efficient usage of an explicit $P_1$ finite-element scheme for the 
solution of coefficient inverse problems (CIPs), in the particular context described above 
was made evident in \cite{BK}.
 In many algorithms aimed at solving electromagnetic CIPs, a qualitative collection of
experimental measurements is necessary on the boundary of a 
computational domain, in order to determine the dielectric permittivity function
therein. In this case, in principle the numerical solution of the time-dependent
Maxwell's equations is required in the entire space
$\mathbb{R}^{3}$ (see e.g. \cite{BK, BCN, BTKM, MalmbergBeilina1,
  MalmbergBeilina2, Malmberg}, but instead it can be more    
efficient to consider Maxwell's equations with a constant dielectric
permittivity in a neighborhood of the boundary of a 
computational domain.  
The explicit $P_1$ finite-element scheme considered in this work was numerically
tested in the solution of the time-dependent Maxwell's system in both two- and three-dimensional 
geometry (cf. \cite{cejm}). It was also combined with a few algorithms to solve different CIPs for
determining the dielectric permittivity function in connection with the time-dependent
Maxwell's equations, using both simulated and experimentally generated data (see \cite{BCN, BTKM, MalmbergBeilina1, MalmbergBeilina2, Malmberg}).
In short, the formal reliability analysis of such a method conducted in this work, corroborates the previously observed adequacy of this numerical approach.\\

An outline of this paper is as follows: In Section 2 we describe in
detail the model problem being solved, and give its equivalent
variational form. In Section 3 we set up the discretizations of the
model problem in both space and time, and recall the main  
results of the reliability analysis conducted in \cite{arxiv1} for the underlying numerical model.
Section 4 is devoted to the numerical experiments that validate such results. We
conclude in Section 5 with a few comments.

\section{The model problem}

The particular form of Maxwell's equations for the electric field ${\bf e}=(e_1,e_2)$ in a bounded domain $\Omega$ of $\Re^2$ with boundary $\partial \Omega$ that we deal with in this work is 
as follows. First we consider that $\Omega = \bar{\Omega}_{in} \cup \Omega_{out}$, where $\Omega_{in}$ is an interior open set whose boundary does not intersect 
$\partial \Omega$ and $\Omega_{out}$ is the complementary set of $\bar{\Omega}_{in}$ with respect to $\Omega$. Now in case ${\bf e}$ satisfies (homogeneous) Dirichlet  boundary conditions, we are given ${\bf e}_0 \in [H^1(\Omega)]^2$ and ${\bf e}_1 \in {\bf H}(div,\Omega)$ satisfying $\nabla \cdot (\varepsilon {\bf e}_0) = \nabla \cdot (\varepsilon {\bf e}_1) = 0$ where $\varepsilon$ is the electric permittivity. $\varepsilon$ is assumed to belong to $W^{2,\infty}(\Omega)$ and to fulfill $\varepsilon \equiv 1 $ in $\Omega_{out}$ and $\varepsilon \geq 1$. Incidentally, throughout this article we denote the standard semi-norm of $C^m(\bar{\Omega})$ by $| \cdot |_{m,\infty}$ for $m >0$ and the standard norm of $C^{0}(\bar{\Omega})$ by $\| \cdot \|_{0,\infty}$. \\
\indent In doing so, the problem to solve is:
\begin{equation}\label{eq1}
  \begin{array}{ll}
    \varepsilon \partial_{tt} {\bf e} + \nabla \times  \nabla \times {\bf e}  = {\bf 0} & \mbox{ in } \Omega \times (0, T), \\
    {\bf e}(\cdot,0) = {\bf e}_0(\cdot), \mbox{ and } \partial_t{\bf e}(\cdot,0) = {\bf e}_1(\cdot) & \mbox{ in } \Omega, \\
    {\bf e} = {\bf 0} & \mbox{ on } \partial \Omega \times (0,T), \\
    \nabla \cdot (\varepsilon {\bf e}) = {\bf 0} & \mbox{ in } \Omega.
  \end{array}
\end{equation}

\begin{remark} As pointed out above the analysis carried out in \cite{arxiv1} applies to the case where absorbing boundary conditions $\partial_n {\bf e}=-\partial_t {\bf e}$ are prescribed, where $\partial_n {\bf e}$ represents the outer normal derivative of ${\bf e}$ on $\partial \Omega $. This choice was motivated by the fact that they correspond to practical situations addressed in \cite{BCN, BTKM, MalmbergBeilina1,MalmbergBeilina2, Malmberg}. 
\rule{2mm}{2mm}
\end{remark} 

We next set \eqref{eq1} in variational form. With this aim we denote the standard inner product of $[L^2(\Omega)]^2$ by $(\cdot,\cdot)$ and the corresponding norm by $\parallel \{\cdot\} \parallel$. Further, for a given non-negative function $\omega \in L^{\infty}(\Omega)$ we introduce the weighted $L^2(\Omega)$-semi-norm 
$\| \{\cdot\} \|_{\omega}:=\sqrt{\int_{\Omega} |\omega| |\{\cdot\}|^2 d{\bf x}}$, which is actually a norm if $\omega \neq 0$ everywhere in $\bar{\Omega}$. We also 
introduce, the notation $({\bf a},{\bf b})_{\omega}:= \int_{\Omega} \omega {\bf a} \cdot {\bf b} d{\bf x}$ for two fields ${\bf a},{\bf b}$ which are square integrable in $\Omega$. Notice that if  
$\omega$ is strictly positive this expression defines an inner product associated with the norm $\| \{\cdot\} \|_{\omega}$.\\
Then requiring that ${\bf e}_{|t=0} = {\bf e}_0$ and $\{\partial_t{\bf e}\}_{|t=0} = {\bf e}_1$ and 
${\bf e}=0$ on $\partial \Omega \times [0,T]$, we write for all $ {\bf v} \in [H_0^1(\Omega)]^2$,
\begin{equation}\label{eq2}
 \left (\partial_{tt} {\bf e},{\bf v} \right )_{\varepsilon} + (\nabla {\bf e},\nabla {\bf v})  + (\nabla \cdot \varepsilon {\bf e}, \nabla \cdot {\bf v}) - 
(\nabla \cdot {\bf e}, \nabla \cdot {\bf v})  = 0 \;\forall t \in (0, T).
\end{equation}

We recall that the equivalence of problem (\ref{eq2}) with Maxwell's equations \eqref{eq1} was established in \cite{arxiv1}. 

\section{The numerical model}

Henceforth we restrict our studies to the case where $\Omega$ is a polygon. 

\subsection{Space semi-discretization} 
Let $V_h$ be the usual $P_1$ FE-space of continuous functions related to a mesh ${\mathcal T}_h$ fitting $\Omega$, consisting of triangles with maximum edge length $h$, belonging to a quasi-uniform family of meshes (cf. \cite{Ciarlet}). \\
Setting ${\bf V}_h := [V_h \cap H^1_0(\Omega)]^2$ we define ${\bf e}_{0h}$ (resp. ${\bf e}_{1h}$) to be the usual ${\bf V}_h$-interpolate of 
${\bf e}_0$ (resp. ${\bf e}_1$). Then the semi-discretized problem is space we wish to solve writes,\\

\emph{Find }${\bf e}_{h} \in {\bf V}_{h}$ \emph{ such  that $\forall {\bf v} \in {\bf V}_{h}$ }
\begin{equation}\label{eq6}
  \begin{array}{l}
    \left ( \partial_{tt} {\bf e}_{h}, {\bf v} \right )_{\varepsilon} + (\nabla {\bf e}_h,\nabla {\bf v})+ (\nabla \cdot[\varepsilon {\bf e}_h], \nabla \cdot {\bf v} ) 
    -(\nabla \cdot {\bf e}_h, \nabla \cdot {\bf v}) = 0,  \\
		\\
   {\bf e}_h(\cdot,0) = {\bf e}_{0h}(\cdot) \mbox{ and } \partial_t {\bf e}_h(\cdot,0) = {\bf e}_{1h}(\cdot) \mbox{ in } \Omega.  
  \end{array}
\end{equation}

\subsection{Full discretization}

To begin with we consider a natural centered time-discretization scheme to solve \eqref{eq6}, namely: Given a number $N$ of time steps we define the time increment $\tau := T/N$. Then we approximate ${\bf e}_h(k\tau)$
by ${\bf e}_h^k \in {\bf V}_h$ for $k=1,2,...,N$ according to the following scheme for $k=1,2,\ldots,N-1$:
\begin{equation}
\label{eq7consist}
  \begin{array}{l}
    \displaystyle \left(\frac{{\bf e}_h^{k+1} - 2 {\bf e}_h^k + {\bf e}_h^{k-1}}{\tau^2}, {\bf v} \right)_{\varepsilon} + (\nabla {\bf e}_h^k, \nabla {\bf v}) 
     + (\nabla \cdot \varepsilon {\bf e}_h^k, \nabla \cdot {\bf v}) - (\nabla \cdot {\bf e}_h^k, \nabla \cdot {\bf v}) = 0 \; \forall {\bf v} \in {\bf V}_h,\\
		\\
     {\bf e}_h^0 = {\bf e}_{0h} \mbox{ and } {\bf e}_h^1 = {\bf e}_h^0 + \tau {\bf e}_{1h} \mbox{ in } \Omega.
  \end{array}
\end{equation}

Owing to its coupling with ${\bf e}_h^{k}$ and ${\bf e}_h^{k-1}$ on the left hand side of \eqref{eq7consist}, ${\bf e}_h^{k+1}$ cannot be determined explicitly by \eqref{eq7consist} at every time step. In order to enable an explicit solution we resort to the classical mass-lumping technique. We recall that for  a constant $\varepsilon$ this 
consists of replacing on the left hand side the inner product $( {\bf u},{\bf v})_{\varepsilon}$ by an inner product 
$({\bf u},{\bf v})_{\varepsilon.h}$, using the trapezoidal rule to compute the integral of $\int_K \varepsilon {\bf u}_{|K} \cdot {\bf v}_{|K}d{\bf x}$ (resp. $\int_{K \cap \partial \Omega} {\bf u}_{|K} \cdot {\bf v}_{|K}dS$), for every element $K$ in ${\mathcal T}_h$, where ${\bf u}$ stands for ${\bf e}_h^{k+1} - 2 {\bf e}_h^k + {\bf e}_h^{k-1}$. It is well-known that in this case the matrix associated with 
$(\varepsilon {\bf e}_h^{k+1},{\bf v})_h$ for ${\bf v} \in {\bf V}_h$, is a diagonal matrix. In our case $\varepsilon$ is not constant, but the same property will hold if we replace in each element $K$ the integral of 
$\varepsilon {\bf u}_{|K} \cdot {\bf v}_{|K}$ in a triangle $K \in {\mathcal T}_h$ as follows:

\[ 
\int_K \varepsilon {\bf u}_{|K} \cdot {\bf v}_{|K} d{\bf x} \approx \varepsilon(G_K)  \displaystyle area(K) \sum_{i=1}^3 \frac{{\bf u}(S_{K,i}) \cdot {\bf v}(S_{K,i})}{3}, 
\]
\noindent where $S_{K,i}$ are the vertexes of $K$, $i=1,2,3$, $G_K$ is the centroid of $K$.\\

Before pursuing we define the auxiliary function $\varepsilon_h$ whose value in each $K \in {\mathcal T}_h$ is constant equal to $\varepsilon(G_K)$. Then still denoting the approximation of ${\bf e}_h(k\tau)$
by ${\bf e}_h^k$, for $k=1,2,...,N$ we determine ${\bf e}_h^{k+1}$ by,
\begin{equation}\label{eq7}
  \begin{array}{l}
    \displaystyle \left( \frac{{\bf e}_h^{k+1} - 2 {\bf e}_h^k + {\bf e}_h^{k-1}}{\tau^2}, {\bf v} \right)_{\varepsilon_h,h} + (\nabla {\bf e}_h^k, \nabla {\bf v}) 
     + (\nabla \cdot \varepsilon {\bf e}_h^k, \nabla \cdot {\bf v}) - (\nabla \cdot {\bf e}_h^k, \nabla \cdot {\bf v}) = 0 \; \forall {\bf v} \in {\bf V}_h,\\
		\\
     {\bf e}_h^0 = {\bf e}_{0h} \mbox{ and } {\bf e}_h^1 = {\bf e}_h^0 + \tau {{\bf e}_1}_h \mbox{ in } \Omega.
  \end{array}
\end{equation}

\subsection{Convergence results}

Recalling the assumption that $\varepsilon \in W^{2,\infty}(\Omega)$ we first set 
\begin{equation}
\label{eta}
\eta :=  2 + | \varepsilon |_{1,\infty} + 2| \varepsilon|_{2,\infty}; \\
\end{equation}
Next we recall the classical inverse inequality (cf. \cite{Ciarlet}) together with a result in \cite{JCAM2010} according to which,
\begin{equation}
\label{inversineq} 
\| \nabla v \| \leq C h^{-1} \| v \|_{\varepsilon_h,h} \mbox{ for all } v \in V_h,
\end{equation} 
where $C$ is a mesh-independent constant. \\
Now we assume that $\tau$ satisfies the following CFL-condition: 
\begin{equation} 
\label{CFL}
\tau \leq h /\nu \mbox{ with } \nu=C(1+ 3 \| \varepsilon - 1 \|_{\infty})^{1/2}.
\end{equation}.
We further assume that the solution ${\bf e}$ to equation \eqref{eq1} belongs to $[H^4\{\Omega \times(0,T)\}]^2$. \\

Let us define a function ${\bf e}_h$ in $\bar{\Omega} \times [0,T]$ whose value at $t=k\tau$ equals ${\bf e}_h^k$ for $k=1,2,\ldots,N$ and that varies linearly with $t$ in each time interval $([k-1] \tau,k \tau)$, in such a way that $\partial_t{\bf e}_{h}({\bf x},t)=\displaystyle \frac{{\bf e}^k_h({\bf x}) - {\bf e}^{k-1}_h({\bf x})}{\tau}$ for every ${\bf x} \in \bar{\Omega}$ and $t \in ([k-1]\tau,k \tau)$. 
We also define ${\bf a}^{m+1/2}(\cdot)$ for any field ${\bf a}(\cdot,t)$ to be ${\bf a}(\cdot,[m+1/2]\tau)$. \\
Then denoting by $| \cdot |_m$ the standard semi-norm of Sobolev space $H^m(\Omega)$ for $m \in {\mathcal N}$, according to \cite{arxiv1} we have:  \\

Provided the CFL condition \eqref{CFL} is fulfilled and $\tau$ also satisfies $\tau \leq 1/[2\eta]$, 
under the above regularity assumption on ${\bf e}$, there exists a constant ${\mathcal C}$ depending only on $\Omega$, $\varepsilon$ and $T$ such that,
  
\begin{equation}
\label{convergence}
\boxed{
\begin{array}{l}
\displaystyle \displaystyle \max_{1 \leq m \leq N-1}  \left\| [\partial_t ({\bf e}_h - {\bf e})]^{m+1/2} \right\| 
+ \displaystyle \displaystyle \max_{2 \leq m \leq N}  \| \nabla ({\bf e}_h^{m} - {\bf e}^{m}) \| \\
\leq {\mathcal C} (\tau + h + h^2/\tau) \displaystyle \left\{ \| {\bf e} \|_{H^4[\Omega \times (0,T)]} + |{\bf e}_0|_2 + |{\bf e}_1 |_2 \right\}. 
\mbox{ \rule{2mm}{2mm}}
\end{array}
}
\end{equation}   
\eqref{convergence} means that, as long as $\tau$ varies linearly with $h$, first order convergence of scheme \eqref{eq7} in terms of either $\tau$ or $h$ holds in the sense of the norms on the left hand side of \eqref{convergence}.

\section{Numerical validation}  

We perform numerical tests in time $(0,T)=(0,0.5)$ in the computational domain $\Omega =
[0,1] \times [0,1]$  for the model problem in two-dimension space, namely  
\begin{equation}\label{model}
  \begin{array}{ll}
    \varepsilon \partial_{tt} {\bf e}  - \nabla^2 {\bf e} 
    - \nabla \nabla \cdot (\varepsilon - 1) {\bf e}   = {\bf f} & \mbox{ in } \Omega \times (0, T), \\
    {\bf e}(\cdot,0) = {\bf 0} \mbox{ and } \partial_t  {\bf e}(\cdot,0)
    = {\bf 0} & \mbox{ in } \Omega, \\
    {\bf e} =  {\bf 0} & \mbox{ on } \partial \Omega \times (0,T).
  \end{array}
\end{equation}
for the electric field  $ {\bf e} = (e_1, e_2)$.

\begin{figure}[h!]
\begin{center}
\begin{tabular}{cccc}
  {\includegraphics[scale=0.1, clip=]{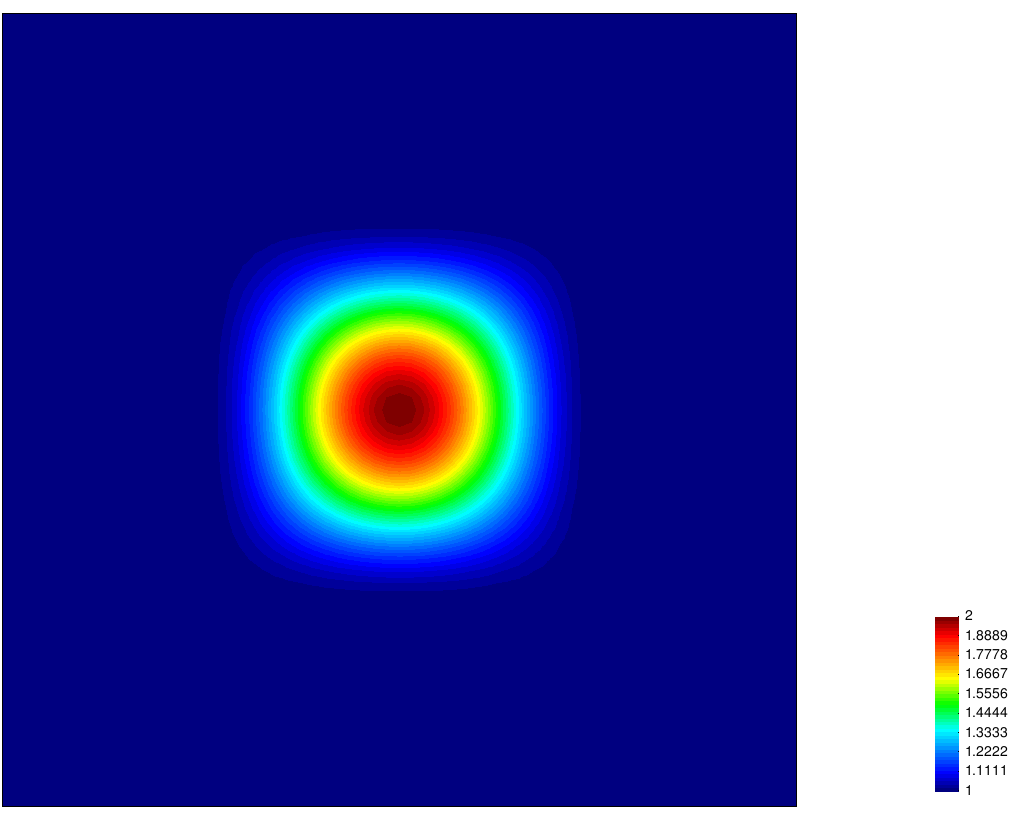}} &
  {\includegraphics[scale=0.1, clip=]{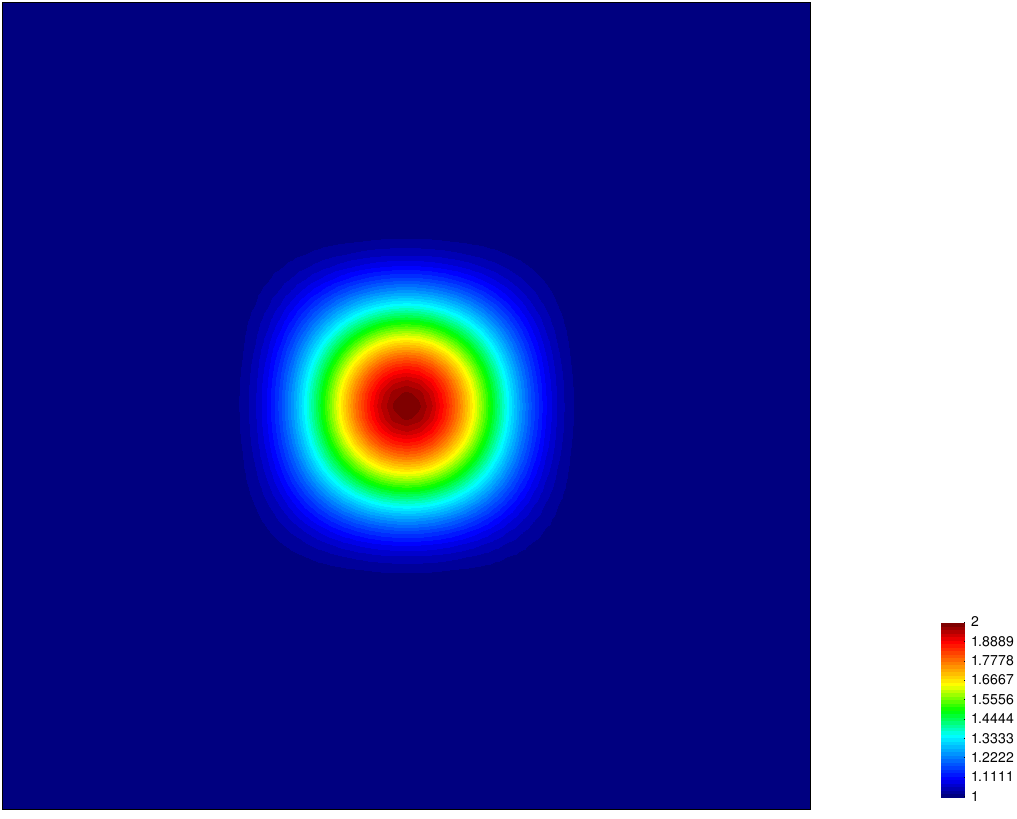}} &
  {\includegraphics[scale=0.1, clip=]{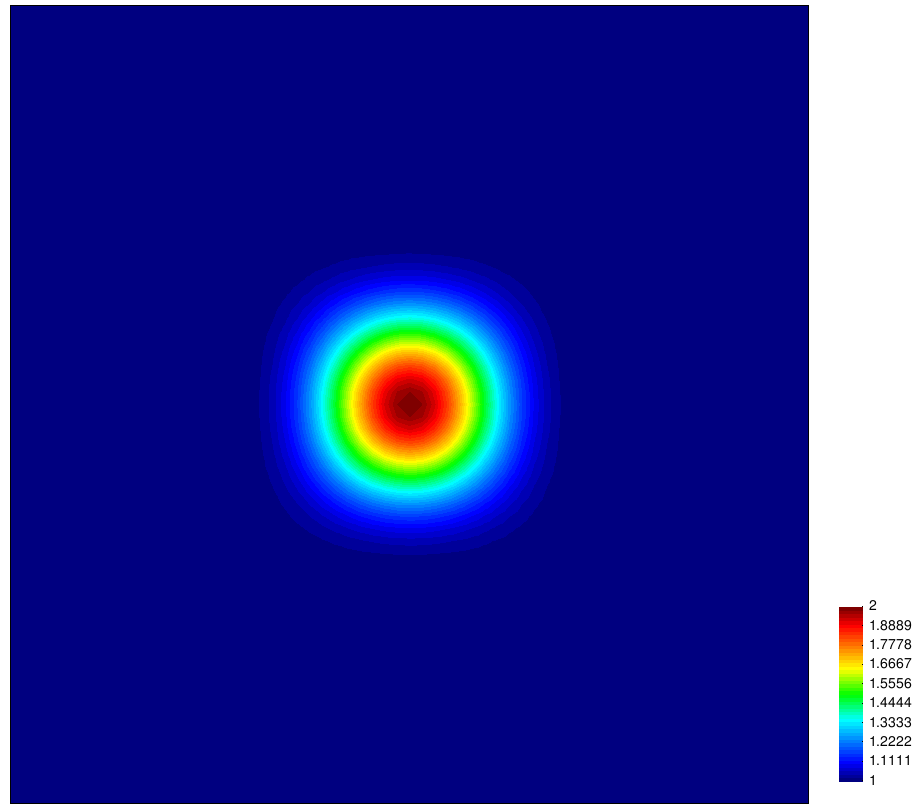}} &
   {\includegraphics[scale=0.1, clip=]{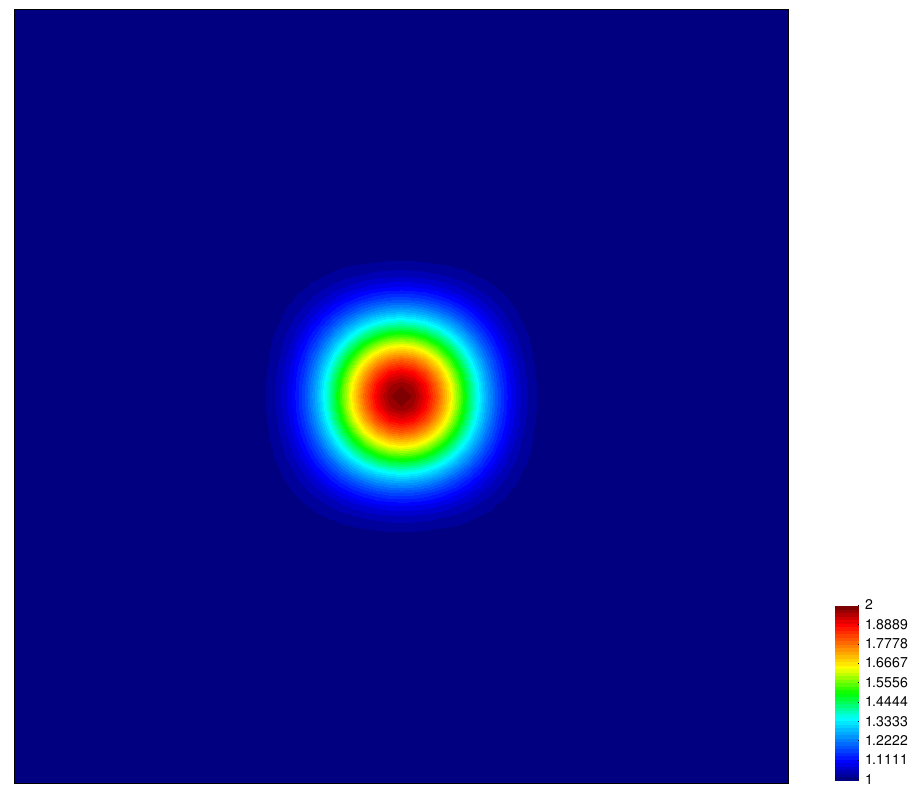}} \\
   $m=2$ & $m=3$ &  $m=4$ &  $m=5$\\
    {\includegraphics[scale=0.1, clip=]{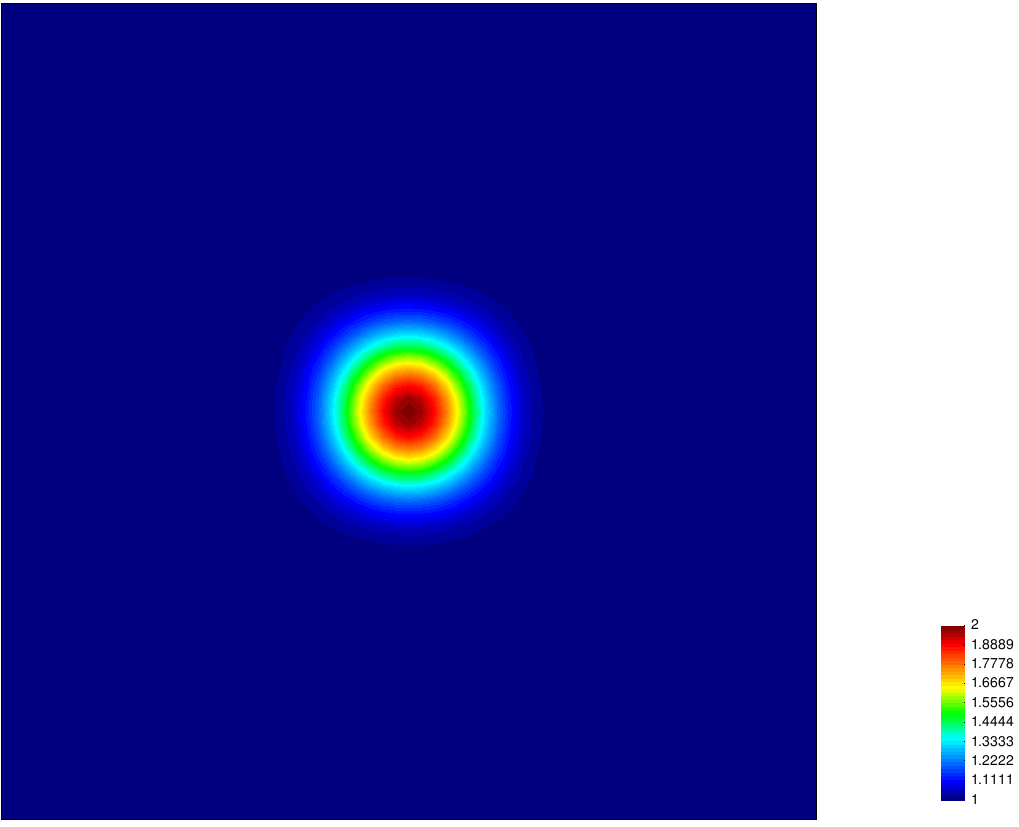}} &
  {\includegraphics[scale=0.1, clip=]{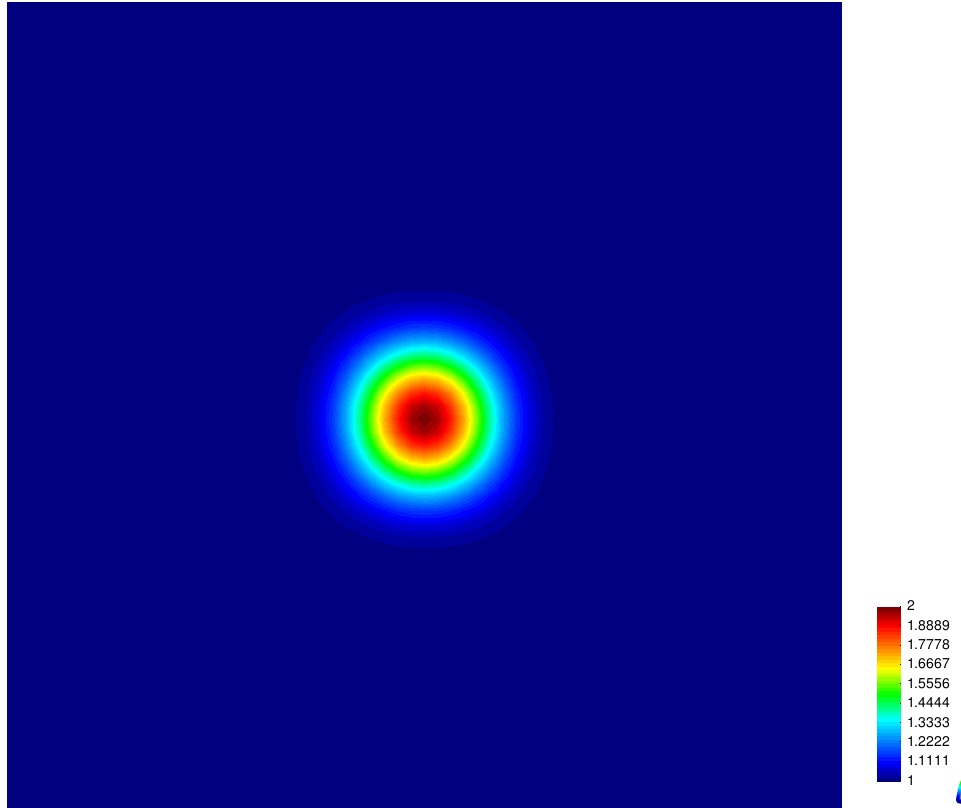}} &
  {\includegraphics[scale=0.1, clip=]{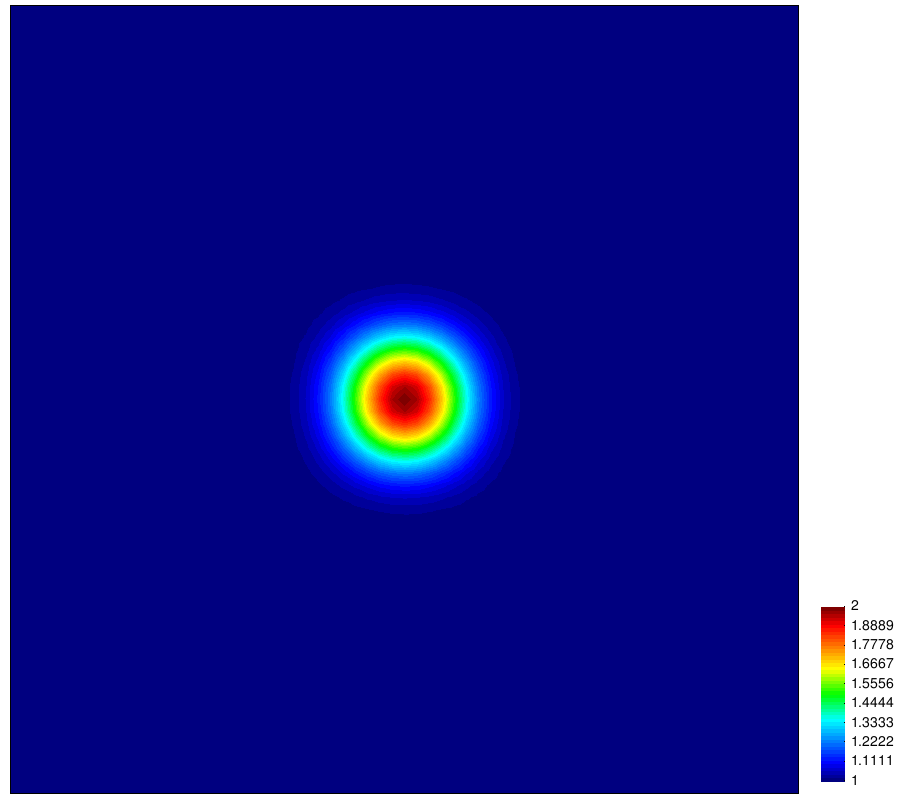}} &
   {\includegraphics[scale=0.1, clip=]{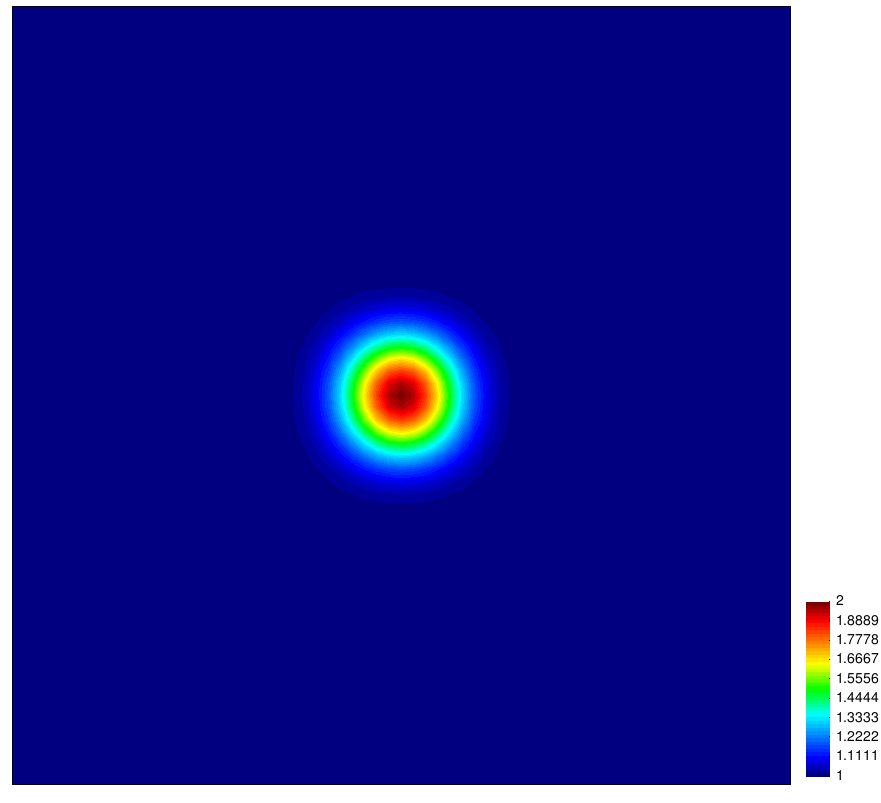}} \\
  $m=6$ & $m=7$ &  $m=8$ &  $m=9$\\
\end{tabular}
\end{center}
\caption{Function $\varepsilon(x,y)$ in the domain $\Omega=[0,1] \times[0,1]$ for different values of $m$ in \eqref{eps}}
\label{fig:F1}
\end{figure}

 The source data ${\bf f}$ (the right hand side) is chosen such that the function
\begin{equation}\label{eq1}
  \begin{split}
    e_1 &=    \frac{1}{\varepsilon} 2 \pi \sin^2 \pi x  \cos \pi y  \sin \pi y ~ \frac{t^2}{2}, \\
   e_2 &=  - \frac{1}{\varepsilon}   2 \pi   \sin^2 \pi y \cos \pi x   \sin \pi x   ~\frac{t^2}{2}
  \end{split}
\end{equation}
is the exact solution of the model problem \eqref{model}
In \eqref{eq1} the function $\varepsilon$ is defined to be,
\begin{equation}\label{eps}
  \varepsilon(x,y)= \left \{
  \begin{array}{ll}
    1 + \sin^m \pi (2x-0.5) \cdot \sin^m \pi (2y-0.5) & \textrm{in $[0.25,0.75] \times [0.25, 0.75]$}, \\
    1 &  \textrm{otherwise},
  \end{array}
  \right.
\end{equation}
where $m$ is an integer greater than one. In Figure \ref{fig:F1} the function $\varepsilon$ are illustrated for different  values of $m$.

The solution given by \eqref{eq1} satisfies homogeneous
initial conditions together with homogeneous Dirichlet conditions on the
boundary $\partial \Omega$ of the square $\Omega$ for every time $t$. In our 
computations we used the software package Waves \cite{waves} only
for the finite element method applied to the solution of the model
problem \eqref{model}. We note that this package was also used in \cite{cejm} to solve the
the same model problem \eqref{model} by a domain
decomposition FEM/FDM method.\\

 We discretized the computational domain $\Omega \times (0,T)$ denoting by
 ${K_h}_l = \{K\}$ a partition of the spatial domain $\Omega$ into triangles $K$ of sizes $h_l= 2^{-l}, l=1,...,6$.
We let $J_{\tau_l}$ be a partition of the time domain $(0,T)$ into time
 intervals $J=(t_{k-1},t_k]$ of uniform length $\tau_l$
for a given number of time intervals $N, l=1,...,6$. We choose the time step
$\tau_l = 0.025 \times 2^{-l}, \; l=1,...,6$, which provides numerical stability for all meshes. 

We performed numerical tests taking $m=2,..., 9$ in \eqref{eps} and computed the 
maximum value over the time steps of the relative errors measured in the function $L_2$-norm and the $H^1$-semi-norm and in the 
$L_2$ norm for the time-derivative, respectively represented by
\begin{equation}\label{norms}
  \begin{split}
e_l^1 &= \displaystyle \frac{\displaystyle \max_{1 \leq k \leq N} \| {\bf e}^k-{\bf e}_h^k \|}{\displaystyle \max_{1 \leq k \leq N}\| {\bf e}^k \|},\\
e_l^2 &= \displaystyle \frac{ \displaystyle \max_{1 \leq k \leq N} \| \nabla ({\bf e}^k - {\bf e}^k_h) \|}{
\displaystyle \max_{1 \leq k \leq N} \| \nabla {\bf e}^k \|}, \\
e_l^3 &=  \displaystyle \frac{\displaystyle \max_{1 \leq k \leq N-1} \| \{ \partial_t ({\bf e}-{\bf e}_h)\}^{k+1/2} \|}{\displaystyle \max_{1 \leq k \leq N-1} \| \{ \partial_t {\bf e}\}^{k+1/2}\|}.
\end{split}
  \end{equation}
Here ${\bf e}$ is the exact solution given by \eqref{eq1} and ${\bf e}_h$  is the computed solution, while $N= T /\tau_l$.\\
In Tables 1-4 method's convergence in these three senses is observed taking $m=2,3,6,7$ in \eqref{eps}.
\begin{table}[h!] 
\center
\begin{tabular}{ l  l l l l l l l l }
\hline
$l$ &  $nel$  & $nno$ &  $e_l^1$  & $e_{l-1}^1/e_{l}^1$ &
$e_l^2$
  & $e_{l-1}^2/e_{l}^2$ & $e_l^3$ &   $e_{l-1}^3/e_{l}^3$\\
\hline 
1 & 8 & 9       &  0.054247   &              & 0.2767 &           & 1.0789  &       \\
2 & 32& 25      &  0.013902   & 3.902100       & 0.1216 &  2.2755   & 0.4811  & 2.2426 \\
3 & 128 & 81    &  0.003706  &   3.751214     & 0.0532 &  2.2857   & 0.2544  &  1.8911 \\
4 & 512 & 289   &  0.000852   & 4.349765      & 0.0234 &  2.2735   & 0.1279  &  1.9891 \\
5 & 2048 & 1089 &  0.000229   & 3.720524       & 0.0121 &  1.9339   & 0.0641  &  1.9953\\
6 & 8192 & 4225 &  0.000059   &  3.881356           & 0.0061 &  1.9836   & 0.0321 &  1.9969 \\
\hline
\end{tabular}
\caption{Maximum over the time steps of relative errors in the $L_2$-norm, in the $H^1$-seminorm and in the $L^2$-norm 
of the time derivative for mesh sizes $h_l= 2^{-l}, l=1,...,6$ taking $m=2$ in \eqref{eps}}
\label{test1}
\end{table}

\begin{table}[h!]
\center
\begin{tabular}{ l  l l l l ll ll }
  \hline
  $l$ &  $nel$  & $nno$ &  $e_l^1$  & $e_{l-1}^1/e_{l}^1$ &
$e_l^2$
  & $e_{l-1}^2/e_{l}^2$ & $e_l^3$ &   $e_{l-1}^3/e_{l}^3$\\
\hline 
1 & 8 & 9       & 0.043394   &           & 0.2784   &           & 1.0869  &       \\
2 & 32& 25      & 0.011451  & 3.789538     & 0.1098   &  2.5355   & 0.5305  & 2.0488 \\
3 & 128 & 81    & 0.003343  & 3.425366    & 0.06     &  1.83     & 0.2586  &  2.0514 \\
4 & 512 & 289   & 0.000781   & 4.385873    & 0.0248   &  2.4194   & 0.1306  &  1.9801 \\
5 & 2048 & 1089 & 0.000202   & 3.866337   & 0.0119   &  2.0840   & 0.0654  &  1.9969\\
6 & 8192 & 4225 & 0.000052   & 3.884615    & 0.0059   &  2.0169   & 0.0327 &  2 \\
\hline
\end{tabular}
\caption{Maximum over the time steps of relative errors in the $L_2$-norm, in the $H^1$-seminorm and in the $L^2$-norm 
of the time derivative for mesh sizes $h_l= 2^{-l}, l=1,...,6$ taking $m=3$ in \eqref{eps}}
\label{test2}
\end{table}

\begin{table}[h!]
\center
\begin{tabular}{ l  l l l l ll ll }
  \hline
  $l$ &  $nel$  & $nno$ &  $e_l^1$  & $e_{l-1}^1/e_{l}^1$ &
$e_l^2$
  & $e_{l-1}^2/e_{l}^2$ & $e_l^3$ &   $e_{l-1}^3/e_{l}^3$\\
\hline 
1 & 8 & 9       & 0.054228   &           & 0.2837   &             & 1.1120  &          \\
2 & 32& 25      & 0.012241  & 4.430030     & 0.0906   &  3.1313     & 0.4937  & 2.2524   \\
3 & 128 & 81    & 0.002973  & 4.117389    & 0.0408    & 2.2206     & 0.2665  &  1.8525  \\
4 & 512 & 289   & 0.000590   & 5.038983    & 0.0150   &  2.7200     & 0.1335  &  1.9963  \\
5 & 2048 & 1089 & 0.000163   & 3.619631    & 0.0079   &  1.8987     & 0.0667  &  2.0015  \\
6 & 8192 & 4225 & 0.000043   & 3.790698   & 0.0040   &  1.9750     & 0.0334 &  1.9970   \\
\hline
\end{tabular}
\caption{Maximum over the time steps of relative errors in the $L_2$-norm, in the $H^1$-seminorm and in the $L^2$-norm 
of the time derivative for mesh sizes $h_l= 2^{-l}, l=1,...,6$ taking $m=6$ in \eqref{eps}}
\label{test3}
\end{table}

\begin{table}[h!]
\center
\begin{tabular}{ l  l l l l ll ll }
  \hline
  $l$ &  $nel$  & $nno$ &  $e_l^1$  & $e_{l-1}^1/e_{l}^1$ &
$e_l^2$
  & $e_{l-1}^2/e_{l}^2$ & $e_l^3$ &   $e_{l-1}^3/e_{l}^3$\\
\hline 
1 & 8 & 9       & 0.054224   &           & 0.5710   &             & 1.1208  &          \\
2 & 32& 25      & 0.012483  & 4.343828     & 0.1505   &  3.7940     & 0.5024 & 2.2309   \\
3 & 128 & 81    & 0.002751  & 4.537623    & 0.0686   &  2.1939     & 0.2688  &  1.8690  \\
4 & 512 & 289   & 0.000627  & 4.387559    & 0.0240   &  2.8583     & 0.1339  &  2.0075  \\
5 & 2048 & 1089 & 0.000158   & 3.968354    & 0.0114   &  2.1053     & 0.0669  &  2.0015  \\
6 & 8192 & 4225 & 0.000040   & 3.949999         & 0.0057   &  2          & 0.0334 &  2.0030   \\
\hline
\end{tabular}
\caption{Maximum over the time steps of relative errors in the $L_2$-norm, in the $H^1$-seminorm and in the $L^2$-norm 
of the time derivative for mesh sizes $h_l= 2^{-l}, l=1,...,6$ taking $m=7$ in \eqref{eps}}
\label{test4}
\end{table}

\begin{figure}[h!]
\begin{center}
\begin{tabular}{cc}
  {\includegraphics[scale=0.5, clip=]{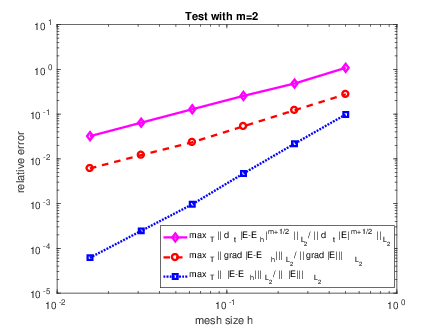}} &
  {\includegraphics[scale=0.5, clip=]{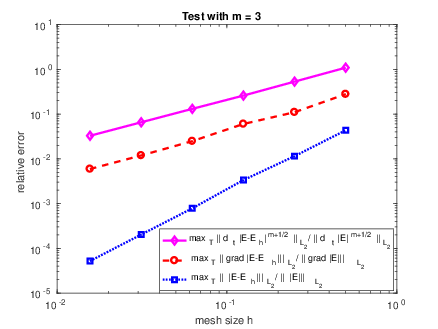}}   
\end{tabular}\
\end{center}
\caption{Maximum in time of relative errors for $m=2$ (left) and $m=3$ (right)}
\label{fig:F2}
\end{figure}

\begin{figure}[h!]
\begin{center}
\begin{tabular}{cc}
  {\includegraphics[scale=0.5, clip=]{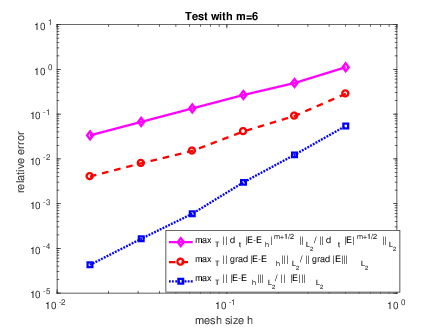}} &
  {\includegraphics[scale=0.5, clip=]{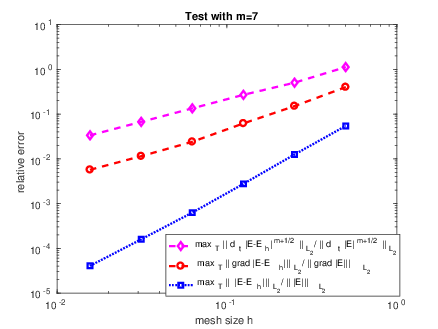}}
\end{tabular}
\end{center}
\caption{Maximum in time of relative errors for $m=6$ (left) and $m=7$ (right)}
\label{fig:F3}
\end{figure}

\begin{figure}[h!]
\begin{center}
\begin{tabular}{cc}
  {\includegraphics[scale=0.5, clip=]{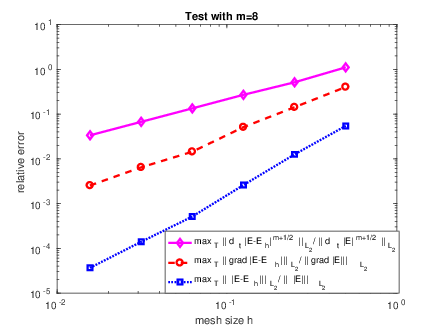}} &
  {\includegraphics[scale=0.5, clip=]{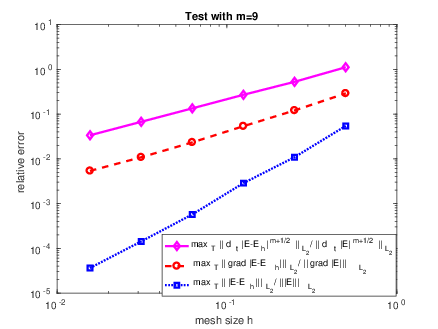}}   
\end{tabular}
\end{center}
\caption{Maximum in time of relative errors for $m=8$ (left) and $m=9$ (right)}
\label{fig:F4}
\end{figure}

Figure \ref{fig:F2} shows convergence of our numerical scheme based on a $P_1$ space discretization, taking the 
function $\varepsilon$ defined by \eqref{eps} with $m=2$ (on the
left) and $m=3$ (on the right) for
$\varepsilon(x)$. Similar convergence results are presented in
Figures \ref{fig:F3} and \ref{fig:F4} taking $m=6,7,8,9$ in \eqref{eps}.

Observation of these tables and figures clearly indicates that our scheme behaves like a first order method in the (semi-)norm of $L^{\infty}[(0,T);H^1(\Omega)]$ for ${\bf e}$ and in the norm of $L^{\infty}[(0,T);L^2(\Omega)]$ for $\partial_t {\bf e}$ for all the chosen values of $m$.  
As far as the values $m=6$ and $m=8$ are concerned this perfectly conforms to the a priori error estimates 
established in \cite{arxiv1}. However those tables and figures also show that such theoretical predictions  
extend to cases not considered in our analysis such as $m=2$ and $m=3$, in which the regularity of the exact solution is lower than assumed, or yet in the cases $m=7$ and $m=9$, in which the minimum of $\varepsilon$ is not attained on the boundary. Otherwise stated some of our assuptions seem to be of academic interest only and a lower regularity of the solution such as $\{H^2[\Omega \times (0,T)]\}^2$ should be sufficient to attain optimal first order convergence in both senses.
On the other hand second-order convergence can be expected from our scheme in the norm of $L^{\infty}[(0,T);L^2(\Omega)]$ for ${\bf e}$, according to Tables 1-4 and Figures 2-4. 

\section{Summary}
In this work we validated the reliability analysis conducted in \cite{arxiv1} for a numerical scheme to solve Maxwell's equations of electromagnetism, combining an explicit finite difference time discretization with a lumped-mass $P1$ finite element space discretization. The scheme is effective in the particular case where the dielectric permittivity is constant in a neighborhood of the boundary of the spacial domain. After presenting the problem under consideration for the electric field we supplied the detailed description of such a scheme and recalled the a priori error estimates that hold for the latter under suitable regularity assumptions specified in \cite{arxiv1}.
Then we showed by means of numerical experiments performed for a test-problem in two-dimension space with known exact solution, that the convergence results given in \cite{arxiv1} are confirmed in practice. Furthermore we presented convincing evidence that such theoretical predictions extend to solutions with much lower regularity than the one assumed in our analysis. Similarly optimal second-order convergence is observed in a norm other than those in which convergence was formally established.   
In short we undoubtedly indicated that Maxwell's equations can be efficiently solved with classical conforming linear finite elements in some relevant particular cases, among which lie the model problem \eqref{model}. \\
   
\newpage

\noindent \textbf{\underline{Acknowledgment:}} The research of the first author
is supported by the Swedish Research Council grant VR 2018-03661. The
second author gratefully acknowledges the financial support provided
by CNPq/Brazil through grant 307996/2008-5. \rule{2mm}{2mm}


\end{document}